\RequirePackage{fix-cm}

\documentclass[article]{IEEEtran}

\usepackage[utf8]{inputenc} 

\usepackage[square,sort&compress,numbers]{natbib}
\usepackage{graphicx}
\usepackage{geometry} 
\usepackage{graphicx} 
\usepackage{mathtools}

\usepackage{xcolor}
\usepackage{booktabs} 
\usepackage{array} 
\usepackage{enumitem} 
\usepackage{verbatim} 
\usepackage{amsmath}
\usepackage{setspace}
\usepackage{times}
\usepackage{multirow}
\usepackage[font=small]{caption}
\usepackage{longtable}
\usepackage{amssymb}

\usepackage{fancyhdr} 
\newenvironment{nospaceflalign}
 {\setlength{\abovedisplayskip}{5pt}\setlength{\belowdisplayskip}{5pt}%
  \csname flalign\endcsname}
 {\csname endflalign\endcsname\ignorespacesafterend}

\usepackage{algorithm}
\usepackage{algpseudocode}
\newcommand{\revision}[1]{{\leavevmode\color{black}#1}}

\begin{document}

\title{Global Solution Strategies for the Network-Constrained Unit Commitment Problem with AC Transmission Constraints}

\author{Jianfeng~Liu,
        Anya~Castillo,~\IEEEmembership{Member,~IEEE,}
        Jean-Paul~Watson,~\IEEEmembership{Member,~IEEE,}
        and~Carl~D.~Laird,~\IEEEmembership{Member,~IEEE}
\thanks{J. Liu is with the Department
of Chemical Engineering, Purdue University, West Lafayette,
IN, 47907 USA e-mail: liu1172@purdue.edu.}
\thanks{A. Castillo is with Sandia National Laboratories, Albuquerque,
NM, 87185 USA e-mail: arcasti@sandia.gov.}
\thanks{J.-P. Watson is with Sandia National Laboratories, Albuquerque,
NM, 87185 USA e-mail: jwatson@sandia.gov.}
\thanks{C.D. Laird is with Davidson School of
  Chemical Engineering, Purdue University, West Lafayette,
IN, 47907 USA and with Sandia National Laboratories, Albuquerque,
NM, 87185 USA e-mail: cdlaird@sandia.gov.}
\thanks{Manuscript received 1/12/18; under review in IEEE TPWRS.}}

\maketitle

\begin{abstract}
We propose a novel global solution algorithm for the network-constrained unit commitment problem incorporating a nonlinear alternating current model of the transmission network, which is a \revision{nonconvex mixed-integer nonlinear programming (MINLP)} problem. Our algorithm is based on the multi-tree \revision{global optimization} methodology, which iterates between a mixed-integer lower-bounding problem and a nonlinear upper-bounding problem. We exploit the mathematical structure of the \revision{unit commitment problem with AC power flow constraints (UC-AC)} and leverage optimization-based bounds tightening, second-order cone relaxations, and piecewise outer approximations to guarantee a globally optimal solution at convergence. Numerical results on four benchmark problems illustrate the effectiveness of our algorithm, both in terms of convergence rate and solution quality. 
\end{abstract}

\subsection{Notation}
\allowdisplaybreaks
\noindent Sets
\begin{align*}
& \mathcal{B} && \text{Set of all buses $\{1, ..., B\}$} \\
& \mathcal{B}_{b} && \text{Set of all buses that are connected to bus $b$} \\
& \revision{\mathcal{C}} && \revision{\text{Set of all cycles in a cycle basis for the network}} \\
& \mathcal{G} && \text{Set of all generators $\{1, ..., G\}$} \\
& \mathcal{G}_{b} && \text{Set of all generators at bus $b$} \\ 
& \mathcal{L} && \text{Set of all branches (transmission lines)}\\
& \mathcal{L}_{c} && \revision{\text{Set of branches in cycle $c$}}\\
& \mathcal{L}_{b}^{in} && \text{Set of all \emph{inbound} branches to bus $b$}  \\
& \mathcal{L}_{b}^{out} && \text{Set of all \emph{outbound} branches from bus $b$} \\
& \mathcal{S}_{g} && \text{Set of startup segments of generator $g$ $\{1, ..., S_{g}\}$} \\
& \mathcal{SC} && \text{Set of all synchronous condensers $\{1, ..., SC\}$} \\
& \mathcal{SC}_{b} && \text{Set of all synchronous condensers at bus $b$} \\ 
& \mathcal{T} && \text{Set of time periods $\{1, ..., T\}$} 
\end{align*}
Parameters
\begin{align*}
& A_{g,n} && \text{Coefficients ($n=0,1,2$) of quadratic}\\
&&& \text{production cost function of generator $g$ }\\
& B_{b}^{sh}&& \text{Shunt susceptance at bus $b$} \\
& B_{l}&& \text{Imag. part of branch $l$ admittance matrix} \\
& G_{b}^{sh}&& \text{Shunt conductance at bus $b$} \\
& G_{l}&& \text{Real part of branch $l$ admittance matrix} \\
& K_{g,\tau}^{su}&& \text{Startup  cost of generator $g$} \\
& K_{g}^{sd} && \text{Shutdown cost of generator $g$} \\
& P_{b,t}^{D} &&\text{Real  power demand at bus $b$, time $t$} \\
& P_{t}^{R} && \text{System reserve requirement at time $t$} \\
& P_{g}^{G,min}&& \text{Min. real power output of generator $g$} \\
& P_{g}^{G,max} && \text{Max. real power output of generator $g$} \\
& Q_{b,t}^{D} &&\text{Reactive power demand at bus $b$, time $t$} \\
& Q_{g}^{G,min}&& \text{Min. reactive power output of generator $g$} \\
& Q_{g}^{G,max} && \text{Max. reactive power output of generator $g$} \\
& Q_{sc}^{SC,min}&& \text{Min. output of synchronous condenser $sc$}\\
& Q_{sc}^{SC,max} && \text{Max. output of synchronous condenser $sc$} \\
& RD_{g}&& \text{Ramp-down limit of generator $g$} \\
&  RU_{g} && \text{Ramp-up limit of generator $g$} \\
& S_{l}^{max} && \text{Apparent power limit on branch $l$} \\
& SD_{g}  && \text{Shutdown capability of generator $g$} \\
& SU_{g} && \text{Startup capability of generator $g$} \\
& T_{g,\tau}^{su} && \text{Startup cost function time segment for} \\
& && \text{generator $g$} \\
& T_{g}^{u} && \text{Min. uptime of generator $g$} \\
& T_{g}^{d} && \text{Min. downtime of generator $g$} \\
& V_{b}^{min}&& \text{Min. voltage magnitude at bus $b$} \\
& V_{b}^{max} && \text{Max. voltage magnitude at bus $b$} 
\end{align*}
Variables
\begin{align*}
& \delta_{g,\tau,t} && \text{Startup cost segment indicator}\\
& \theta_{l,t}&& \text{Voltage phase angle difference between ends}\\
&&&\text{(bus $b$ and bus $k$) of branch $l$ at time $t$, $\theta_{b,k,t}$} \\
& c_{b,k,t}&& \text{Second-order cone variable}\\
& c_{g,t}^p && \text{Production cost for generator $g$ at time $t$}  \\
& f^p && \text{Total production cost} \\
&  f^{sd} && \text{Total shutdown cost} \\
&  f^{su} && \text{Total startup cost} \\
& p_{g,t}^{G}&& \text{Real power output of generator $g$ at time $t$} \\
& p_{l,t}^{f}&& \text{Real power flow \emph{from} branch $l$, at time $t$}\\
& p_{l,t}^{t} &&\text{Real power flow \emph{to} branch $l$, at time $t$}\\
& q_{g,t}^{G} && \text{Reactive power output of generator $g$ at time $t$} \\
& q_{l,t}^{f} && \text{Reactive power flow \emph{from} branch $l$, at time $t$}\\
& q_{l,t}^{t} && \text{Reactive power flow \emph{to} branch $l$, at time $t$}\\
& q_{sc,t}^{SC} && \text{Reactive power output of synchronous}\\
&&&\text{condenser $sc$ at time $t$} \\
& r_{g,t}^{a} && \text{Real power reserve provided by generator $g$}\\
&&&\text{at time $t$} \\
& s_{b,k,t} && \text{Second-order cone variable}\\
& u_{g,t} && \text{Startup status, equal to 1 if generator} \\
&&&\text{$g$ starts up at time $t$, 0 otherwise} \\
& v_{b,t} && \text{Voltage magnitude at bus $b$ at time $t$,}\\
&&&\text{$v_{b,t}^2 = (v_{b,t}^{r})^2 +  (v_{b,t}^{j})^2$} \\
& v_{b,t}^{j} && \text{Imag. part of voltage phasor at bus $b$, time $t$} \\
& v_{b,t}^{r}&& \text{Real part of voltage phasor at bus $b$, time $t$} \\
& w_{g,t} && \text{Shutdown status, equal to 1 if generator} \\
&&&\text{$g$ shuts down at time $t$, 0 otherwise} \\
& y_{g,t} && \text{Unit on/off status, equal to 1 if generator}\\
&&&\text{$g$ is on-line at time $t$, 0 otherwise} 
\end{align*}

\section{Introduction}
\label{sec:introduction}
\IEEEPARstart{R}{ecently} \revision{the Federal Energy Regulatory Commission (FERC) reported that uplift, which represents out-of-market payments that result when an out-of-merit generation cost is incurred to relieve a constraint, can arise due to the inability of independent system operators (ISOs) to fully model the physical constraints on an alternating current (AC) network \cite{Sauer2014}.  Recent work on the day-ahead unit commitment problem, which was led by MISO technical staff in \cite{chen2016}, attests to the importance and non-trivial complexity of incorporating network constraints due to the performance challenges introduced by denser matrices and additional nonlinearities.} 

\revision{Because of these modeling difficulties, current practice is to perform unit commitment using DC approximations (or copper plate) to represent the transmission network. These approximations do not allow rigorous treatment of AC power flow constraints. As a result,
  certain resources are consistently committed outside of the market to address unforeseen reliability issues; this results in concentrated uplift payments \cite{Sauer2014}. Such resources are often required for reactive power compensation in order to provide system voltage control that enables more efficient delivery and utilization of real power \cite{Ilic2015}. Because such reliability requirements are largely unmodeled in day-ahead unit commitment, more cost effective resources are displaced for these out-of-merit commitments. Alternatively, in the real-time market, operators may have to manually commit and dispatch reliability units while also manually re-dispatching or de-committing other resources, e.g., exceptional dispatches in CAISO \cite{CAISO}, out-of-merit generation in NYISO \cite{NYISO}, and balancing operating reserves in PJM \cite{PJM}.} 

\revision{To address these concerns, this paper focuses on solution of the unit commitment problem with AC power flow constraints (UC-AC). Solving real-world operations and market settlement with alternating current optimal power flow (ACOPF) is not trivial. Due to the scale of real-world power systems, network-constrained unit commitment problems can be extremely large and computationally challenging to solve. Coupling this with nonconvex AC powerflow constraints leads to a mixed-integer nonlinear programming (MINLP) problem that is NP-hard \cite{Verma2009,Lehmann2015,Tseng1996}. If the continuous relaxation of the MINLP is a convex optimization problem, we refer to it as a \emph{convex} MINLP. Otherwise, the problems is referred to as a \emph{nonconvex} MINLP. With this definition, the UC-AC is a nonconvex MINLP. Algorithms exist to address both convex and nonconvex MINLP problems, however, tailored solution strategies are often required to achieve desired computational performance. In this paper, we present the first known global optimization approach that can successfully solve the UC-AC on a set of small- to medium-sized test problems.}


\revision{Deterministic MINLP algorithms can be classified into single-tree and multi-tree methods.} Single-tree deterministic algorithms, i.e., the well-known branch-and-bound (BB) methods~\cite{Land1960, Dakin1965}, seek a global optimum by searching a single tree using a systematic enumeration strategy consisting of three primary steps: branching, bounding, and selection. BB-based global optimization strategies have been well-studied and specialized, yielding strategies such as Branch-and-Reduce \cite{Ryoo1996}, Reduced Space Branch-and-Bound~\cite{Epperly1997}, Branch-and-Contract~\cite{Zamora1999}, Branch-and-Cut \cite{Kesavan2000}, and Branch-and-Sandwich~\cite{Kleniati2014}. These approaches are suitable for general, nonconvex MINLP problems of small or medium size, but become computationally intractable with increasing numbers of discrete variables (like those arising in UC-AC). 

In contrast, multi-tree methods~\cite{Belotti2013} iteratively solve a sequence of related lower-bounding \revision{(master)} and upper-bounding problems. For \revision{convex} MINLP problems, many multi-tree solution strategies -- including Generalized Benders Decomposition (GBD)~\cite{Geoffrion1972}, Outer Approximation (OA)~\cite{Duran1986,Fletcher1994}, and Exact Cutting Plane (ECP) methods~\cite{Westerlund1995} -- are effective, and have been applied to a broad range of MINLPs in various application domains.
\revision{For \textit{nonconvex} MINLP problems with special properties (e.g., those that are bilinear, polynomial, linear fractional, or concave separable), extensions of these basic multi-tree methods have been reported in the literature~\cite{Porn1999,Porn2000}. Despite their importance, Bonami, K{\i}l{\i}n\c{c} and Linderoth noted that recent advancements in respective MILP and NLP problem classes have unfortunately resulted in ``far more modest" improvements in general algorithms for even convex MINLPs \cite{Bonami2012}, illustrating the need for specialized approaches.}

\revision{The classic OA approach, a multi-tree technique, was originally developed to solve \emph{convex} MINLP. This approach solves a sequence of MILP master and convex NLP subproblems and yields a globally optimal solution for a convex MINLP in a finite number of iterations for a given $\epsilon$-tolerance on the optimality gap~\cite{Duran1986,Fletcher1994}. The MILP master problem is a relaxation of the original MINLP that provides a provable lower bound on the MINLP along with a candidate integer solution. Fixing the integers in the MINLP yields a convex NLP subproblem that provides a valid upper bound and a candidate solution (for both continuous and integer variables) to the overall MINLP. In this classic approach, the master problem is further refined (i.e., relaxation strengthened) though the addition of linear outer approximations of convex constraints in the MINLP. The algorithm iterates between the master problem and the NLP subproblem, and terminates when the gap between the lower and upper bounds is sufficiently closed. Constraints can also be added to the master problem to remove previously visited integer solutions (using so-called \emph{integer cuts}). These methods have been extended to nonconvex problems where global convergence of the MINLP can be achieved as long as global solutions of the NLP subproblem are ensured \cite{ignacio-paper-from-jianfeng}.  Similar multi-tree solution strategies for nonconvex MINLP has also been successfully used in various applications \cite{Bergamini2005, Bergamini2007, Karuppiah2008}.}

\revision{We extend \cite{Liu2017a} and propose a multi-tree method based on OA for the UC-AC problem. The master problem is constructed using second-order cone (SOC) relaxations of the nonconvex AC transmission constraints \cite{Jabr2007}. As the algorithm iterates, the master problem is further refined with piecewise outer approximations to strengthen the tightness of the relaxation and the \emph{lower-bound} computation. The algorithm from \cite{Liu2017a} is used to find a global solution of the nonconvex NLP in the \emph{upper-bound} computation. Furthermore, we incorporate optimization-based bounds tightening (OBBT) techniques that are valid in both master and subproblem iterations and, since our proposed approach enforces global solution of the NLP subproblem, we are able to include \emph{integer cuts} in the master problem that remove previously visited solutions from the feasible space as the algorithm iterates.} To the best of our knowledge, this is the first global solution algorithm \revision{successfully applied to the UC-AC problem, identifying} solutions with quality certificates (optimality gaps) in time-limited environments.

The remainder of this paper is organized as follows. In Section~\ref{sec:ncuc} we introduce the unit commitment formulation with AC transmission constraints (UC-AC). In Section~\ref{sec:multi-tree} we outline the necessary problem relaxations and our global optimization algorithm.  In Section~\ref{sec:results} we report numerical results on a variety of test systems. We then conclude in Section~\ref{sec:conclusions} with a summary of our contributions and directions for future work.

\section{UC-AC Problem Formulation}
\label{sec:ncuc}
We now introduce our UC-AC problem formulation. We first present the core UC model in Section \ref{sec:uc}, which is based on the compact three-binary (3BIN) formulation introduced in~\cite{Morales2013a}. We then present the rectangular power-voltage (RPQV) model \cite{Cain2013} in Section \ref{sec:acopf} to represent the steady-state operations of the nonlinear AC transmission network. We integrate these constraint sets to represent the UC-AC problem, resulting in a \revision{nonconvex MINLP}. A tailored solution technique for this model is proposed in the following section.

\subsection{Unit Commitment Model}\label{sec:uc}
We use the term \emph{UC skeleton} when referring to a unit commitment model consisting only of a cost function, operating constraints, and any associated continuous and binary variables with no network representation. We summarize several key components of the \emph{3BIN} formulation here; refer to \cite{Morales2013a} for further details.

\subsubsection{Cost Function} 

The total cost in UC is the sum of three major components -- production costs, startup costs, and shutdown costs -- as follows:
\begin{equation*}
f^p + f^{su} + f^{sd}.
\end{equation*}
We assume that the production cost $f^{p}$ is a quadratic monotonically non-decreasing function of real power generation; in practice, this is often replaced with a piecewise approximation. Computation of $f^{p}$ in the quadratic case is accomplished by imposing the constraints
\begin{nospaceflalign}\label{eq:Total_Production_Cost_1}
 & A_{g,2}(p_{g,t}^{G})^2 + A_{g,1} p_{g,t}^{G} + A_{g,0} y_{g,t} \le c_{g,t}^{p} && \forall \; g, \; t 
\end{nospaceflalign}
\begin{nospaceflalign}\label{eq:Total_Production_Cost_2}
 & f^{p} = \sum\nolimits_{g \in \mathcal{G}} \sum\nolimits_{t \in \mathcal{T}}c_{g,t}^{p}&&
\end{nospaceflalign}
where $A_{g,2}$, $A_{g,1}$, and $A_{g,0}$ are known cost coefficients in ($\$/$MW$^2$h), ($\$/$MWh) and ($\$/$h)  associated with a specific generator $g$.

To formulate the total startup cost, $f^{su}$, we first introduce a new binary variable $\delta_{g,\tau,t}$, which indicates the startup type $\tau$ of generator $g$ at time period $t$. In particular, $\delta_{g,\tau,t}$ takes the value of $1$ if the generator $g$ starts up at time $t$ and has been previously offline within $[T_{g,\tau}^{su}, T_{g,\tau+1}^{su})$ hours. The logical constraints between $w_{g,t}$, $u_{g,t}$, and $\delta_{g,\tau,t}$ are given as
\begin{nospaceflalign}\label{eq:Startup_Cost_1}
 & \delta_{g,\tau.t} \le \sum\nolimits_{t'=t-T_{g,\tau}^{su}}^{t+1-T_{g,\tau+1}^{su}}w_{g,t'} && \forall \; g, \; t, \; \tau \in [1, S_{g}) 
\end{nospaceflalign}
\begin{nospaceflalign}\label{eq:Startup_Cost_2}
 & u_{g,t} = \sum\nolimits_{\tau \in \mathcal{S}_{g}} \delta_{g,\tau,t} && \forall \; g, \; t
\end{nospaceflalign}
where $S_{g}$ is the number of startup types for generator $g$, and $u_{g,t}$ and $w_{g,t}$ indicate startup and shutdown of generator $g$ in time $t$, respectively. Note that $w_{g,t}$ with positive time index $t$ are variables, otherwise $w_{g,t}$ are treated as constants to demonstrate previous system status. 

For a thermal unit, the startup cost is assumed to be a monotonically increasing step function with respect to the generator's previous off-line time. The total startup cost is given by
\begin{nospaceflalign}\label{eq:Total_Startup_Cost}
 & f^{su} = \sum\nolimits_{g \in \mathcal{G}} \sum\nolimits_{t \in \mathcal{T}} \sum\nolimits_{\tau \in \mathcal{S}_{g}} K_{g,\tau}^{su} \delta_{g,\tau,t}&&
\end{nospaceflalign}
where $K_{g,\tau}^{su}$ is the cost of startup type $\tau$ for generator $g$. Given logical constraints \eqref{eq:Startup_Cost_1} and \eqref{eq:Startup_Cost_2}, and the monotonically non-decreasing startup cost function, it can be shown that $\delta_{g,\tau,t}$ will always solve to a binary value. In other words, instead of explicitly defining $\delta_{g,\tau,t}$ as a binary, it can be relaxed as a continuous variable within range $[0,1]$.

The shutdown cost of generator $g$ is assumed to be independent of its previous on-line states, and the total shutdown cost is:
\begin{nospaceflalign} \label{eq:Total_Shutdown_Cost}
 & f^{sd} = \sum\nolimits_{g \in \mathcal{G}} \sum\nolimits_{t \in \mathcal{T}} K_{g}^{sd} w_{g,t}.&&
\end{nospaceflalign}

\subsubsection{Operating Constraints} 

According to operating restrictions, a thermal unit must stay in one state (either on-line or off-line) for a certain period of time before its state can be changed again. Such time periods vary between different generator types. To enforce this requirement, we have to introduce minimum uptime and downtime constraints
\begin{nospaceflalign}\label{eq:Minimum_Up/Downtime_1}
& \sum\nolimits_{t' = t - T_{g}^{u} + 1}^{t} u_{g,t'} \le y_{g,t} && \forall \; g, \; t 
\end{nospaceflalign}
\begin{nospaceflalign}\label{eq:Minimum_Up/Downtime_2}
& \sum\nolimits_{t' = t - T_{g}^{d} + 1}^{t} w_{g,t'} \le 1 - y_{g,t} && \forall \; g, \; t 
\end{nospaceflalign}
where $u_{g,t}$ and $w_{g,t}$ with positive time index $t$ are unknown variables, otherwise they are treated as constants to indicate previous system status. Additional constraints are required to denote the logical correlation between $u_{g,t}$, $w_{g,t}$, and $y_{g,t}$ in
\begin{nospaceflalign}\label{eq:Logical}
 & y_{g,t} - y_{g,t-1} = u_{g,t} - w_{g,t} && \forall \; g, \; t .
\end{nospaceflalign}
Note that these constraints ensure that a generator cannot start up and shut down within the same time period. Given the fact that $y_{g,t}$ is a binary variable, imposing constraints \eqref{eq:Minimum_Up/Downtime_1}, \eqref{eq:Minimum_Up/Downtime_2} and \eqref{eq:Logical} together guarantees that $u_{g,t}$ and $w_{g,t}$ take binary values only. Consequently, $u_{g,t}$, $w_{g,t}$, and $\delta_{g,\tau, t}$, though initially defined as binaries, can be relaxed as continuous within $[0,1]$, leaving the $y_{g,t}$ as the only binary variables in our \emph{UC skeleton} formulation. 
The spinning reserve constraint is defined as 
\begin{nospaceflalign}\label{eq:Reserve}
& P_{t}^{R} \le \sum\nolimits_{g \in \mathcal{G}} r_{g,t} && \forall \; t 
\end{nospaceflalign}
and determines the extra generating capacity available by generators included in the commitment solution at time $t$; typically, the spinning reserve is defined as a fraction of the current total power demand. The upper- and lower-bounds of generator output is dependent on its operating state; the real power productions are constrained by $[P_{g}^{G,min}, P_{g}^{G,max}]$, the startup and shutdown capabilities $SD_{g}$ and $SU_{g}$, and state indicators $y_{g,t}$, $u_{g,t}$, and $w_{g,t}$ where both real power generation $p_{g,t}$ and spinning reserve $r_{g,t}$ are accounted for in
\begin{nospaceflalign}
& p_{g,t} + r_{g,t} \le (P_{g}^{G,max}-P_{g}^{G,min})y_{g,t}&&\nonumber
\end{nospaceflalign}
\begin{nospaceflalign}\label{eq:Power_1}
&\quad - (P_{g}^{G,max}-SU_{g})u_{g,t}  &&\forall \; g, \; t
\end{nospaceflalign}
\begin{nospaceflalign}
& p_{g,t} + r_{g,t} \le (P_{g}^{G,max}-P_{g}^{G,min})y_{g,t} &&\nonumber
\end{nospaceflalign}
\begin{nospaceflalign}\label{eq:Power_2}
&\quad- (P_{g}^{G,max}-SD_{g})w_{g,t+1}  &&\forall \; g, \; t
\end{nospaceflalign}
when $T_{g}^{u} = 1$, and 
\begin{nospaceflalign}
& p_{g,t} + r_{g,t} \le (P_{g}^{G,max}-P_{g}^{G,min})y_{g,t}- (P_{g}^{G,max}&&\nonumber
\end{nospaceflalign}
\begin{nospaceflalign}\label{eq:Power_3}
&\quad -SU_{g})u_{g,t}- (P_{g}^{G,max}-SD_{g})w_{g,t+1} &&\forall \; g, \; t
\end{nospaceflalign}
when $T_{g}^{u} \ge 2$. The real power production is also constrained by ramp-up and ramp-down limits, which are given as
\begin{nospaceflalign}\label{eq:Ramping_1}
& p_{g,t} + r_{g,t} - p_{g,t-1} \le RU_{g} &&  \forall \; g, \; t 
\end{nospaceflalign}
\begin{nospaceflalign}\label{eq:Ramping_2}
& -p_{g,t} + p_{g,t-1} \le RD_{g} && \forall \; g, \; t.
\end{nospaceflalign}
Then, the reactive power productions are only constrained by $[Q_{g}^{G,min}, Q_{g}^{G,max}]$ and $y_{g,t}$ in 
\begin{nospaceflalign}\label{eq:Power_4}
& Q_{g}^{G,min}y_{g,t} \le q_{g,t}^{G} \le Q_{g}^{G,max}y_{g,t} && \forall \; g, \; t. 
\end{nospaceflalign}
Synchronous condensers are not modeled with startup/shutdown costs  and their reactive power output is constrained by $[Q_{sc}^{SC,min}, Q_{sc}^{SC,max}]$
\begin{nospaceflalign}\label{eq:Power_5}
& Q_{sc}^{SC,min} \le q_{sc,t}^{SC} \le Q_{sc}^{SC,max} && \forall \; sc, \; t. 
\end{nospaceflalign}

\subsection{AC Transmission Network Model}\label{sec:acopf}
In electric power system analysis, the \emph{RPQV} model is widely-used to represent an AC transmission network; this approach explicitly models real and reactive power flows in terms of complex voltages in the rectangular form.  A transmission line is denoted as $l {\equiv} (b,k)$, where $b$ is the index of the bus at the \emph{from} end and $k$ is the index of the bus at the \emph{to} end of branch $l$. For integration into our \emph{UC skeleton}, the \emph{RPQV} model is given by
\begin{nospaceflalign}\label{eq:M-PQV_1}
& \sum\nolimits_{l \in \mathcal{L}_{b}^{in}}p_{l,t}^{t} + \sum\nolimits_{l \in \mathcal{L}_{b}^{out}}p_{l,t}^{f} + G_{b}^{sh}v_{b,t}^2  \nonumber\\
&\quad+P_{b,t}^{D}- \sum\nolimits_{g \in \mathcal{G}_{b}}p_{g,t}^{G} = 0 & \forall \; b, \; t 	
\end{nospaceflalign}
\begin{nospaceflalign}
& \sum\nolimits_{l \in \mathcal{L}_{b}^{in}}q_{l,t}^{t} + \sum\nolimits_{l \in \mathcal{L}_{b}^{out}}q_{l,t}^{f} - B_{b}^{sh}v_{b,t}^2+Q_{b,t}^{D}&& \nonumber
\end{nospaceflalign}
\begin{nospaceflalign}\label{eq:M-PQV_2}
&\quad - \sum\nolimits_{g \in \mathcal{G}_{b}}q_{g,t}^{G} - \sum\nolimits_{sc \in \mathcal{SC}_{b}}q_{sc,t}^{SC} = 0 & \forall \; b, \; t 
\end{nospaceflalign}
\begin{nospaceflalign}\label{eq:M-PQV_3}
& p_{l,t}^{f} = G_{l}^{ff}v_{b,t}^2 + G_{l}^{ft}(v_{b,t}^{r}v_{k,t}^{r}+v_{b,t}^{j}v_{k,t}^{j}) \nonumber\\
&\quad- B_{l}^{ft}(v_{b,t}^{r}v_{k,t}^{j}-v_{b,t}^{j}v_{k,t}^{r}) & \forall \; l, \; t 
\end{nospaceflalign}
\begin{nospaceflalign}\label{eq:M-PQV_4}
& q_{l,t}^{f} = -B_{l}^{ff}v_{b,t}^2 - B_{l}^{ft}(v_{b,t}^{r}v_{k,t}^{r}+v_{b,t}^{j}v_{k,t}^{j}) \nonumber\\
&\quad- G_{l}^{ft}(v_{b,t}^{r}v_{k,t}^{j}-v_{b,t}^{j}v_{k,t}^{r}) & \forall \; l, \; t 
\end{nospaceflalign}
\begin{nospaceflalign}\label{eq:M-PQV_5}
& p_{l,t}^{t} = G_{l}^{tt}v_{k,t}^2 + G_{l}^{tf}(v_{k,t}^{r}v_{b,t}^{r}+v_{k,t}^{j}v_{b,t}^{j}) \nonumber\\
&\quad- B_{l}^{tf}(v_{k,t}^{r}v_{b,t}^{j}-v_{k,t}^{j}v_{b,t}^{r}) & \forall \; l, \; t 
\end{nospaceflalign}
\begin{nospaceflalign}\label{eq:M-PQV_6}
& q_{l,t}^{t} = -B_{l}^{tt}v_{k,t}^2 - B_{l}^{tf}(v_{k,t}^{r}v_{b,t}^{r}+v_{k,t}^{j}v_{b,t}^{j}) \nonumber\\
&\quad- G_{l}^{tf}(v_{k,t}^{r}v_{b,t}^{j}-v_{k,t}^{j}v_{b,t}^{r}) & \forall \; l, \; t 
\end{nospaceflalign}
\begin{nospaceflalign}\label{eq:M-PQV_7}
& (V_{b}^{min})^2 \le v_{b,t}^2\le (V_{b}^{max})^2 & \forall \; b, \; t 
\end{nospaceflalign}
\begin{nospaceflalign}\label{eq:M-PQV_8}
&  (p_{l,t}^{f})^2 + (q_{l,t}^{f})^2 \le (S_{l}^{max})^2 & \forall \; l, \; t 
\end{nospaceflalign}
\begin{nospaceflalign}\label{eq:M-PQV_9}
& (p_{l,t}^{t})^2 + (q_{l,t}^{t})^2 \le (S_{l}^{max})^2 & \forall \; l, \; t 
\end{nospaceflalign}
where $v_{b,t}^2 {\equiv} (v_{b,t}^{r})^2 + (v_{b,t}^{j})^2$; see \cite{Liu2017a} for details on computing $G_{l}$ and $B_{l}$ branch admittance submatrices. Note that the \emph{RPQV} problem is nonconvex due to bilinear terms and nonconvex quadratics.
 
\subsection{UC-AC Problem Formulation}\label{sec:ncuc_subsection}
\revision{The UC-AC is a nonconvex MINLP formulations that combines the UC skeleton with the nonlinear ACOPF constraints, giving:}
\begin{equation}\label{eq:NCUC}
\begin{aligned}
  &&& \min  f^{p} + f^{su} + f^{sd}& \\  
&&& \text{s.t.} & \\	
&&& \eqref{eq:Total_Production_Cost_1}-\eqref{eq:M-PQV_9}& \\	
&&& \revision{  y_{g,t},u_{g,t},w_{g,t}\in\{0,1\}}&\revision{\forall \; g, \; t } \\
\end{aligned}
\end{equation} 
In the next section we exploit the special mathematical structure of this problem to solve the problem globally.

\section{UC-AC Global Solution Framework}
\label{sec:multi-tree}

\revision{The UC-AC is a nonconvex MINLP, and our proposed algorithm is a \emph{nested} multi-tree method where both the outer and inner algorithm are based on a nonconvex OA approach that solves a sequence of lower-bounding master problems and upper-bounding subproblems. In this section, we first provide a high-level explanation of the nested multi-tree approach used to solve the UC-AC MINLP problem, followed by a detailed description of the master and NLP subproblems and the algorithm definition. Here, we denote $d{=}\left[y, u, w\right]$ to represent the discrete decisions (i.e., generator commitment variables), and $x$ to represent the continuous variables in the UC-AC problem.}
\subsection{Overview}
\revision{Figure~\ref{fig.ucacalg} shows the multi-tree approach for the UC-AC problem. The algorithm iterates between a master problem and an NLP subproblem, and each pair of such solves comprise a major iteration $q$ for candidate solution denoted as $\left[d^q, x^q\right]$.
\begin{figure}[!ht]
  \begin{center}
    \includegraphics[width=2in]{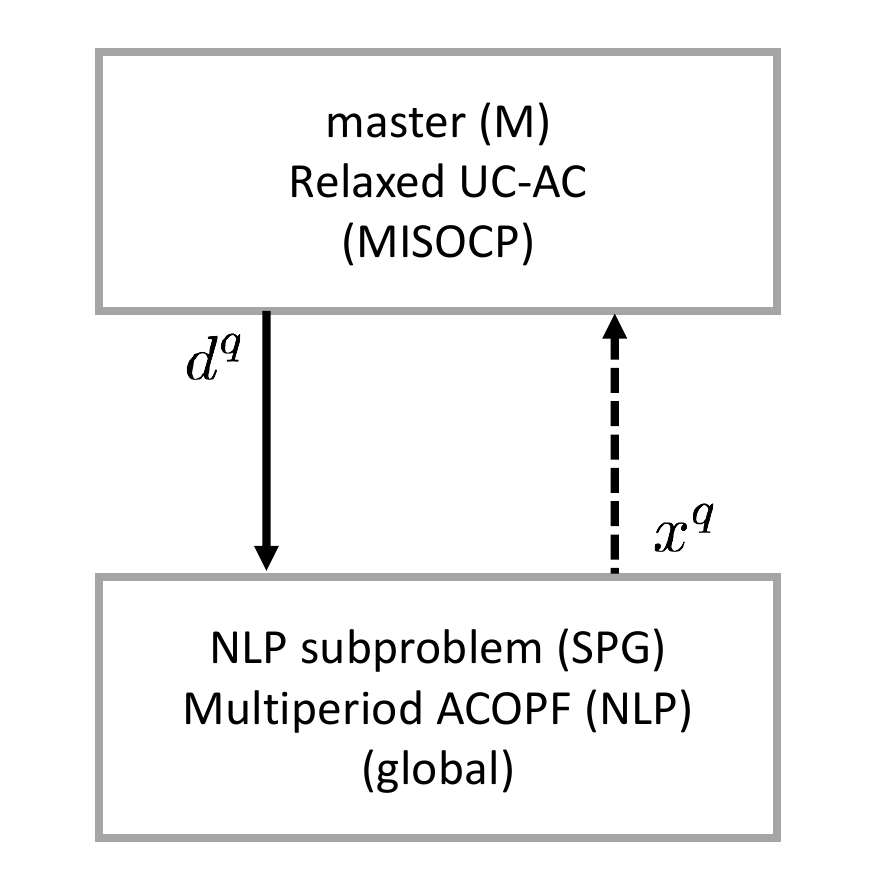}
    \caption{High-level description of the multi-tree approach for global solution of the UC-AC MINLP problem.}
  \label{fig.ucacalg}
  \end{center}
  \end{figure}
The high-level description of the \emph{Outer Algorithm} is as follows: 

The master problem $\mathsf{(M)}$ is a relaxation of the UC-AC problem where the AC power flow constraints are relaxed using the SOC representation from \cite{Jabr2007}. The initial solution of $\mathsf{(M)}$ provides a lower bound on the UC-AC problem and a candidate solution for the binary variables (the generator commitments) given by $d^q$ for iteration $q$. Fixing these variables in the UC-AC MINLP problem yields a nonconvex NLP that represents a multi-period ACOPF problem given by $\mathsf{(SPG)}$. This NLP subproblem, if feasible, provides an upper bound, $z_U^q$, and a candidate solution to the UC-AC, $\left[d^q, x^q\right]$. If the gap between the upper and lower bound is sufficiently small, then the solution has been found, i.e., $z^*{=}z_U^q$ for $\left[d^*, x^*\right]{=}\left[d^q, x^q\right]$.

To further accelerate exploration of the generator commitments, it is also desirable to add cuts to $\mathsf{(M)}$ that remove previously visited solutions $d^q$ from the feasible space. With these \emph{integer cuts} (see Section \ref{sec:integer_cuts}), the solution $z^q_L$ of $\mathsf{(M)}$ is not a true lower bound to the original MINLP, and to ensure convergence with this approach, it is required that we find a globally optimal solution to the NLP subproblem $\mathsf{(SPG)}$ for each candidate binary solution $d^q$. Note that, in the limit, this will result in full enumeration, ensuring convergence of the discrete decision space in a finite number of iterations. However, for the applications and test cases presented in this work, only a few outer iterations were required to close the gap.

}

  \begin{figure}[!ht]
    \begin{center}
    \includegraphics[width=2in]{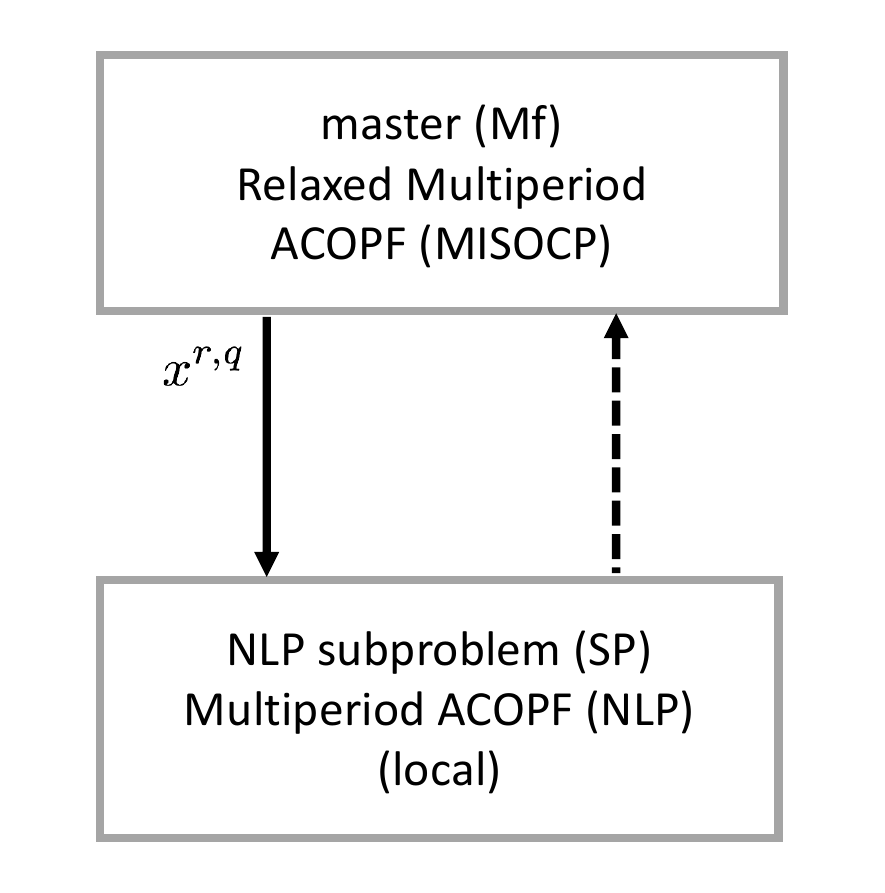}
    \caption{High-level description of the multi-tree approach for global solution of the NLP subproblem $\mathsf{(SPG)}$.}
    \label{fig.spalg}
    \end{center}
  \end{figure}
  
\revision{For global solution of the multi-period ACOPF in $\mathsf{(SPG)}$ we apply the approach of \cite{Liu2017a}, and for completeness, Figure \ref{fig.spalg} shows this algorithm. This strategy is also a multi-tree approach, and hence we refer to the overall algorithm as a \emph{nested} multi-tree approach. Recall that the candidate generator commitments $d^q$ are fixed for this problem. Similar to the \emph{Outer Algorithm} in Figure \ref{fig.ucacalg}, this approach iterates between the master and the NLP subproblem, and each pair of such solves constitutes a minor iteration $r$ on iteration $q$. The high-level description of the \emph{Inner Algorithm} is as follows: 

The master problem $\mathsf{(Mf)}$ is a MISOCP relaxation of the problem $\mathsf{(SPG)}$ ($d^q$ fixed). Therefore, in $\mathsf{(Mf)}$ the only binary variables are those corresponding to piecewise outer approximations. The master problem $\mathsf{(Mf)}$ is solved to find a lower bound for $\mathsf{(SPG)}$, and the solution $x^{r,q}$ from $\mathsf{(Mf)}$ is used to initialize the NLP subproblem $\mathsf{(SP)}$. This NLP subproblem, if feasible provides an upper bound, $z^{r,q}_U$, and a candidate solution $x^{r,q}$. Note that the NLP subproblem $\mathsf{(SP)}$ in Figure \ref{fig.spalg} is the same formulation as $\mathsf{(SPG)}$ in Figure \ref{fig.ucacalg}, however, in this case we only seek a local solution of the NLP subproblem $\mathsf{(SP)}$.

Since we do not add integer cuts to the master problem $\mathsf{(Mf)}$, it is a true relaxation of $\mathsf{(SPG)}$, and closure of the gap between the upper and lower bounds is sufficient to indicate convergence. At each iteration $r$, the master problem is progressively refined by the addition and/or tightening of piecewise outer approximations, as well as optimization-based bounds tightening (OBBT), as discussed later in Sections~\ref{sec:refinements}.

Note that for both \emph{Outer} and \emph{Inner Algorithms}, the respective master problems $\mathsf{(M)}$ and $\mathsf{(Mf)}$ can be further refined with any selection of piecewise outer approximations (see Sections \ref{sec:reverse_cone} and \ref{sec:cycle_constraint}) and with domain reduction techniques, e.g., OBBT (see Section \ref{sec:obbt}). 
}

\revision{
\subsection{Problem Formulations}
\label{sec:formulations}
This section provides a description of the problem formulations $\mathsf{(M)}$, $\mathsf{(SPG)}$,  $\mathsf{(SP)}$, and $\mathsf{(Mf)}$ used in the global algorithm. The master problem $\mathsf{(M)}$ for the UC-AC problem is based on the SOC relaxation of the power flow equations from \cite{Jabr2007}. We replace the quadratic and bilinear terms
in \eqref{eq:NCUC} for all $l{\equiv}(b,k)$ and $t$ with
\begin{nospaceflalign}\nonumber
 c_{b,b,t} &\coloneqq  (v_{b,t}^{r})^2+(v_{b,t}^{j})^2  
\end{nospaceflalign}
\begin{nospaceflalign}\nonumber
 c_{b,k,t} &\coloneqq  v_{b,t}^{r}v_{k,t}^{r}+v_{b,t}^{j}v_{k,t}^{j} 
\end{nospaceflalign}
\begin{nospaceflalign}\nonumber
 s_{b,k,t} &\coloneqq  v_{b,t}^{r}v_{k,t}^{j}-v_{k,t}^{r}v_{b,t}^{j} 
 \end{nospaceflalign}
and introduce a second-order cone relaxation of the condition
\begin{nospaceflalign}\label{eq:SOCR_EQUIV_1}
c_{b,k,t}^2 + s_{b,k,t}^2= c_{b,b,t}c_{k,k,t} &&
 \end{nospaceflalign}
 as
\begin{nospaceflalign}\label{eq:SOCR_EQUIV_1a}
c_{b,k,t}^2 + s_{b,k,t}^2\leq c_{b,b,t}c_{k,k,t}. &&
\end{nospaceflalign}

\subsubsection{\underline{Master Problem $\mathsf{(M)}$}}
With the definitions above, the problem formulation for $\mathsf{(M)}$ is given as follows:
\begin{subequations}
\renewcommand{\theequation}{$\mathsf{M}$.\arabic{equation}}
\label{MISOCP}
\begin{nospaceflalign}\label{eq:RPQV-R_obj}
& z_L \coloneqq \min  f^{p} + f^{su} + f^{sd} &&
\end{nospaceflalign}
\begin{nospaceflalign}
\text{s.t.} \nonumber&&
\end{nospaceflalign}
\begin{nospaceflalign}\label{eq:RPQV-R_0}
\eqref{eq:Total_Production_Cost_1}-\eqref{eq:Power_5},\eqref{eq:M-PQV_8},\eqref{eq:M-PQV_9}&&
\end{nospaceflalign}
\begin{nospaceflalign}\nonumber
\sum\nolimits_{l \in \mathcal{L}_{b}^{in}}p_{l,t}^{t} + \sum\nolimits_{l \in \mathcal{L}_{b}^{out}}p_{l,t}^{f} + G_{b}^{sh}c_{b,b,t}  &&
\end{nospaceflalign}
\begin{nospaceflalign}\label{eq:RPQV-R_1}
\quad+P_{b,t}^{D}- \sum\nolimits_{g \in \mathcal{G}_{b}}p_{g,t}^{G} = 0 &\quad\forall \; b, \; t &
\end{nospaceflalign}
\begin{nospaceflalign}\nonumber
\sum\nolimits_{l \in \mathcal{L}_{b}^{in}}q_{l,t}^{t} + \sum\nolimits_{l \in \mathcal{L}_{b}^{out}}q_{l,t}^{f} - B_{b}^{sh}c_{b,b,t} +Q_{b,t}^{D}&&
\end{nospaceflalign}
\begin{nospaceflalign}\label{eq:RPQV-R_2}
\quad - \sum\nolimits_{g \in \mathcal{G}_{b}}q_{g,t}^{G} - \sum\nolimits_{sc \in \mathcal{SC}_{b}}q_{sc,t}^{SC} = 0& \quad\forall \; b, \; t &
\end{nospaceflalign}
\begin{nospaceflalign}\label{eq:RPQV-R_3}
p_{l,t}^{f} = G_{l}^{ff}c_{b,b,t}+ G_{l}^{ft}c_{b,k,t} - B_{l}^{ft}s_{b,k,t} &\quad \forall \; l, \; t &
\end{nospaceflalign}
\begin{nospaceflalign}\label{eq:RPQV-R_4}
q_{l,t}^{f} = -B_{l}^{ff}c_{b,b,t} - B_{l}^{ft}c_{b,k,t} - G_{l}^{ft}s_{b,k,t}& \quad\forall \; l, \; t&
\end{nospaceflalign}
\begin{nospaceflalign}\label{eq:RPQV-R_5}
p_{l,t}^{t} = G_{l}^{tt}c_{k,k,t}  + G_{l}^{tf}c_{k,b,t} - B_{l}^{tf}s_{k,b,t}&\quad \forall \; l, \; t &
\end{nospaceflalign}
\begin{nospaceflalign}\label{eq:RPQV-R_6}
q_{l,t}^{t} = -B_{l}^{tt}c_{k,k,t}  - B_{l}^{tf}c_{k,b,t} - G_{l}^{tf}s_{k,b,t} &\quad \forall \; l, \; t &
\end{nospaceflalign}
\begin{nospaceflalign}\label{eq:RPQV-R_7}
(V_{b}^{min})^2 \le c_{b,b,t}\le (V_{b}^{max})^2  &\quad\forall \; b, \; t &
\end{nospaceflalign}
\begin{nospaceflalign}\label{eq:RPQV-R_8}
  c_{b,k,t} = c_{k,b,t}&\quad  \forall \; l, \; t  &
\end{nospaceflalign}
\begin{nospaceflalign}\label{eq:RPQV-R_9}
s_{b,k,t} = -s_{k,b,t} & \quad\forall \; l, \; t  &
\end{nospaceflalign}
\begin{nospaceflalign}\label{eq:RPQV-R_10}
c_{b,k,t}^2 + s_{b,k,t}^2 \le c_{b,b,t}c_{k,k,t} &\quad \forall \; l, \; t&
\end{nospaceflalign}
\begin{nospaceflalign}\label{eq:RPQV-R_11}
  y_{g,t},u_{g,t},w_{g,t}\in\{0,1\}&\quad  \forall \; g, \; t&
\end{nospaceflalign}
\end{subequations}

\subsubsection{\underline{NLP Subproblems $\mathsf{(SPG)}$ and $\mathsf{(SP)}$}}
The same NLP subproblem is used in both the outer and the inner multi-tree algorithms, however, for $\mathsf{(SPG)}$, a
global solution is required. The NLP subproblem is
formed by fixing the binary variables $d{=}[y,u,w]$ (generator commitments) in the original
MINLP formulations for the UC-AC. This produces a multi-period ACOPF formulation. For any iteration $j$,
problem for fixed $d^{(j)}{=}[y^{(j)},u^{(j)},w^{(j)}]$ is given as:
\begin{subequations}
\renewcommand{\theequation}{$\mathsf{SP}$}
 \label{eq:NLP}
\begin{nospaceflalign}\nonumber
z_{U} \coloneqq \min  f^{p} + f^{su} + f^{sd} & 
\end{nospaceflalign} 
\begin{nospaceflalign}\nonumber
\text{s.t.} & &
\end{nospaceflalign} 
\begin{nospaceflalign}
\eqref{eq:Total_Production_Cost_1}-\eqref{eq:M-PQV_9}&
\end{nospaceflalign} 
\begin{nospaceflalign}
\text{where} & &\nonumber
\end{nospaceflalign} 
\begin{nospaceflalign}\nonumber
y_{g,t}\coloneqq y^{(j)}_{g,t},u_{g,t}\coloneqq u^{(j)}_{g,t},w_{g,t}\coloneqq w^{(j)}_{g,t} \quad\forall \; g, \; t  &
\end{nospaceflalign} 
\end{subequations} 
\\

\begin{algorithm}[h]
\caption{\revision{\emph{Outer Algorithm} for UC-AC}}
\label{ncuc-algorithm}
\begin{algorithmic}[1]
\State \textbf{Initialization.} 
\newline \quad Iteration $q{=}0$,
\newline \quad $z^*_L \leftarrow -\infty.\;  z^*_U \leftarrow +\infty.\; (d^*, x^*) \leftarrow \O$.
\State \textbf{Solve the Master Problem $(\mathsf{M})$.}\label{marker}
\newline  \quad Solve problem $(\mathsf{M})$ to compute its objective value $z^q_L$ and binary solution $d^q$.
\newline \quad (a) If $(\mathsf{M})$ is infeasible, then $(d^*, x^*)$ is the optimal solution (unless $(d^*, x^*)\equiv\O$, then the UC-AC problem is infeasible). Terminate.
\newline \quad (b) If $z^*_L > z^q_L$, then $z^*_L \leftarrow z^q_L$. 
\State \textbf{Solve for the Upper-Bound.}
\newline \quad Solve the NLP subproblem $(\mathsf{SPG})$ (with fixed $d^q$) to global optimality using \textbf{Algorithm \ref{nested-algorithm}}. Let $z^q_U$ and $(d^q, x^q)$ be the optimal objective value and solution.
\newline \quad (a) If feasible and $z^*_U < z^q_U$, then update the candidate solution: $z^*_U \leftarrow z^q_U$ and $(d^*,x^*) \leftarrow (d^q,x^q)$. 
\State \textbf{Convergence Check}
\newline\quad  (a) If gap $(z^*_U - z^*_L)/z^*_L < \epsilon_O$, the optimal solution $(d^*,x^*)$ has been identified. Terminate. 
\newline\quad  (b) Otherwise add an \emph{integer cut} $(\mathsf{IC})$ for $d^q$ to $(\mathsf{M})$.
\State  \textbf{Iterate} $q \leftarrow q+1$. \textbf{Go to} \textbf{Step \ref{marker}}. 
\end{algorithmic}
\end{algorithm}

\subsubsection{\underline{Master Problem $\mathsf{(Mf)}$}}
Problem $\mathsf{(Mf)}$ is the master problem used in the inner multi-tree approach for obtaining globally optimal solutions to the NLP subproblem $\mathsf{(SPG)}$ from the outer problem. It is based on the same SOC relaxation that is used for problem $\mathsf{(M)}$, however, the generator commitments  $d{=}[y,u,w]$ are fixed. Problem $\mathsf{(Mf)}$ for any iteration $j$ with fixed $d^{(j)}{=}[y^{(j)},u^{(j)},w^{(j)}]$ is given by:
\begin{subequations}
\renewcommand{\theequation}{$\mathsf{Mf}$}
 \label{eq:SOCP}
\begin{nospaceflalign}\nonumber
z_{L_{fixed}}\coloneqq \min  f^{p} + f^{su} + f^{sd} &
\end{nospaceflalign} 
\begin{nospaceflalign}\nonumber
\text{s.t.} & &
\end{nospaceflalign} 
\begin{nospaceflalign}
\eqref{eq:RPQV-R_0}-\eqref{eq:RPQV-R_10} & 
\end{nospaceflalign} 
\begin{nospaceflalign}
\text{where} & &\nonumber
\end{nospaceflalign} 
\begin{nospaceflalign}\nonumber
y_{g,t}\coloneqq y^{(j)}_{g,t},u_{g,t}\coloneqq u^{(j)}_{g,t},w_{g,t}\coloneqq w^{(j)}_{g,t}  \quad\forall \; g, \; t  &
\end{nospaceflalign} 
\end{subequations} 
}

\begin{algorithm}[h!]
\caption{\revision{\emph{Inner Algorithm} for $(\mathsf{SPG})$}}
\label{nested-algorithm}
\begin{algorithmic}[1]
\State \textbf{Initialization.} 
\newline \quad For outer iteration $q$ and fixed binary $d^q$:
\newline \quad Inner iteration $r=0$. 
\newline \quad $z^*_{L_{fixed}} \leftarrow -\infty.\;  z^q_U \leftarrow +\infty.\; x^{q,r} \leftarrow \O$. 
\State \textbf{Solve for the Lower-Bound.}\label{marker-sub}
\newline \quad Solve problem $(\mathsf{Mf})$ (with fixed $d^q$) to find lower bound $z^r_{L_{fixed}}$ solution $x^{q,r}$.
\newline \quad (a) If $(\mathsf{Mf})$ is infeasible then the subproblem $(\mathsf{SPG})$ is infeasible. Return to \textbf{Step 3} in \textbf{Algorithm \ref{ncuc-algorithm}}.
\newline \quad (b) If $z^*_{L_{fixed}} > z^r_{L_{fixed}}$, then $z^*_{L_{fixed}} \leftarrow z^r_{L_{fixed}}$.
\State \textbf{Solve for the Upper-Bound.}
\newline \quad Solve problem $(\mathsf{SP})$ (initialized from $x^{q,r}$) to compute its objective value $z^r_{U_{fixed}}$ and solution $x^{q,r}_{fixed}$. If $z^q_U < z^r_{U_{fixed}}$, then $z^q_U \leftarrow z^r_{U_{fixed}}$ and $x^q \leftarrow x^{q,r}_{fixed}$.
\State \textbf{Convergence Check.}
\newline\quad (a) If $(z^q_U - z^*_{L_{fixed}})/z^*_{L_{fixed}}< \epsilon_I$ (optimality tolerance), then $x^q$ is optimal. Return $z^q_U$ and $x^q$ to \textbf{Step 3} in \textbf{Algorithm \ref{ncuc-algorithm}}.  
\newline \quad (b) Else perform OBBT on selected variables and add or refine partitions for piecewise outer relaxations $(\mathsf{UE})$, $(\mathsf{OE})$, and $(\mathsf{CC})$.
\State \textbf{Iterate} $r \leftarrow r+1$. \textbf{Go to} \textbf{Step \ref{marker-sub}}.
\end{algorithmic}
\end{algorithm}

\revision{
\subsection{Global Solution Algorithm}
\label{sec:multi-tree-algorithm}
In this section, we formally present the \emph{nested} multi-tree algorithm. Algorithm 1 presents the \emph{Outer Algorithm} for the solution of the UC-AC problem, and Algorithm 2 presents the \emph{Inner Algorithm} for global solution of the NLP subproblem from the \emph{Outer Algorithm}. For implementation details on the \emph{integer cuts} $(\mathsf{IC})$, piecewise outer relaxations $(\mathsf{UE})$, $(\mathsf{OE})$, and $(\mathsf{CC})$, and OBBT referred to in the presented algorithms, please see the following section.
}

 \revision{
\subsection{Algorithm Details}
\label{sec:refinements}
\subsubsection{Integer Cuts}\label{sec:integer_cuts}
At each iteration $q$ of the \emph{Outer Algorithm} we add integer cuts that remove previously visited solutions $d^q$. These cuts are given by,
\begin{nospaceflalign}\label{eq:integer-cut}\tag{$\mathsf{IC}$}
&   \sum_{(g,t)\in \mathcal{B}^{(q)}} y_{g,t} - \sum_{(g,t)\in \mathcal{N}^{(q)}} y_{g,t}  \leq   \vert \mathcal{B}^{(q)} \vert - 1 &&
\end{nospaceflalign}
for $q = 1\ldots Q-1$ where $\mathcal{B}^{(q)} =\{g,t \vert y^{(q)}_{g,t} = 1 \}$ and $\mathcal{N}^{(q)} =\{g,t \vert y^{(q)}_{g,t} = 0 \}$ This enhancement ensures that distinct solutions are obtained during each major iteration $q$ of our global solution algorithm.

\subsubsection{``Reverse Cone"}\label{sec:reverse_cone}
For any solution of $(\mathsf{Mf})$, we may have that equation \ref{eq:SOCR_EQUIV_1} is violated, i.e.,
\begin{nospaceflalign}\nonumber
c_{b,b,t}c_{k,k,t} - (c_{b,k,t}^2 + s_{b,k,t}^2) > \varepsilon
\end{nospaceflalign} 
for any $l$ and $t$ due to the second-order cone relaxation of \eqref{eq:SOCR_EQUIV_1}. Therefore, we introduce piecewise relaxations of
\begin{nospaceflalign}\label{eq:SOCR_EQUIV_1b}
c_{b,k,t}^2 + s_{b,k,t}^2 \geq c_{b,b,t}c_{k,k,t}, 
\end{nospaceflalign}
as necessary in each iteration of the \emph{Inner Algorithm}. To describe these relaxations, we define new variables 
\begin{nospaceflalign}\nonumber
 cs_{b,k,t} &\coloneqq  c_{b,k,t}^2 + s_{b,k,t}^2
\end{nospaceflalign}
\begin{nospaceflalign}\nonumber
 cc_{b,k,t} &\coloneqq  c_{b,b,t}c_{k,k,t}
\end{nospaceflalign}
where we construct piecewise over-estimators for $c_{b,k,t}^2 + s_{b,k,t}^2$ and piecewise under-estimators for $c_{b,b,t}c_{k,k,t}$ to obtain an adjustable approximation of \eqref{eq:SOCR_EQUIV_1b}.

Specifically, as first introduced in \cite{Liu2017a}, we extend the bivariate partitioning scheme in \cite{Hasan2010}. We denote our partitioning variables as $cs_{b,k,t,}^{i,j}$ and $cc_{b,k,t}^{i,j}$, where $[\underline{c}_{b,k,t}^{i}, \overline{c}_{b,k,t}^{i}]$ refers to the $i$-th interval for $c_{b,k,t}\in[\underline{c}_{b,k,t}, \overline{c}_{b,k,t}]$ and $[\underline{s}_{b,k,t}^{j}, \overline{s}_{b,k,t}^{j}]$ refers to the $j$-th interval for $s_{b,k,t}\in[\underline{s}_{b,k,t}, \overline{s}_{b,k,t}]$.

The piecewise over-estimators for $cs_{b,k,t}$ are
\begin{subequations}
\renewcommand{\theequation}{$\mathsf{OE}$}
 \label{eq:OE}
\begin{nospaceflalign}\nonumber
cs_{b,k,t}^{i,j} \leq (\underline{c}_{b,k,t}^{i} + \overline{c}_{b,k,t}^{i})c_{b,k,t}^{i,j} + (\underline{s}_{b,k,t}^{j} + \overline{s}_{b,k,t}^{j})s_{b,k,t}^{i,j} &&
\end{nospaceflalign} 
\begin{nospaceflalign}\nonumber
\quad - (\underline{c}_{b,k,t}^{i}\overline{c}_{b,k,t}^{i})c_{b,k,t}^{i,j} + \underline{s}_{b,k,t}^{j}\overline{s}_{b,k,t}^{j})\sigma_{b,k,t}^{i,j} \quad \forall \; (i,j),\; l,\;t&&
\end{nospaceflalign} 
\begin{nospaceflalign}\nonumber
cs_{b,k,t} = \sum\nolimits_{(i,j) \in \Omega^{cs}_{b,k,t}} cs_{b,k,t}^{i,j}  \quad \forall \; l,\; t &&
\end{nospaceflalign} 
\begin{nospaceflalign}\nonumber
\underline{c}_{b,k,t}^{i} \sigma_{b,k,t}^{i,j} \leq c_{b,k,t}^{i,j} \leq \overline{c}_{b,k,t}^{i} \sigma_{b,k,t}^{i,j} \quad \forall \; (i,j),\; l,\; t &&
\end{nospaceflalign} 
\begin{nospaceflalign}\nonumber
c_{b,k,t} = \sum\nolimits_{(i,j) \in \Omega^{cs}_{b,k}} c_{b,k,t}^{i,j} \quad \forall \; l,\; t  &&
\end{nospaceflalign} 
\begin{nospaceflalign}
\underline{s}_{b,k,t}^{j} \sigma_{b,k,t}^{i,j} \leq s_{b,k,t}^{i,j} \leq \overline{s}_{b,k,t}^{j} \sigma_{b,k,t}^{i,j} \quad \forall \; (i,j),\; l,\; t &&
\end{nospaceflalign} 
\begin{nospaceflalign}\nonumber
s_{b,k,t} = \sum\nolimits_{(i,j) \in \Omega^{cs}_{b,k}} s_{b,k,t}^{i,j} \quad \forall \; l,\; t  &&
\end{nospaceflalign} 
\begin{nospaceflalign}\nonumber
\sum\nolimits_{(i,j) \in \Omega^{cs}_{b,k}} \sigma_{b,k,t}^{i,j} =  1 \quad \forall \; l,\; t  &&
\end{nospaceflalign} 
\begin{nospaceflalign}\nonumber
 \sigma_{b,k,t}^{i,j} \in \{0,1\} \quad \forall \; (i,j),\; l,\; t  &&
\end{nospaceflalign} 
\end{subequations}
where $(i,j)\in\Omega^{cs}_{b,k,t} \coloneqq [\underline{c}_{b,k,t}^{i}, \overline{c}_{b,k,t}^{i}] \times [\underline{s}_{b,k,t}^{j}, \overline{s}_{b,k,t}^{j}]$. Then, the piecewise under-estimators for $cc_{b,k,t}$ are
\begin{subequations}
\renewcommand{\theequation}{$\mathsf{UE}$}
 \label{eq:UE}
\begin{nospaceflalign}\nonumber
cc_{b,k,t}^{i,j} \leq \overline{cc}_{k,k,t}^{i,j}cc_{b,b,t}^{i,j} + \overline{cc}_{b,b,t}^{i,j}cc_{k,k,t}^{i,j} &&
\end{nospaceflalign} 
\begin{nospaceflalign}\nonumber
\quad - \overline{cc}_{b,b,t}^{i,j} \overline{cc}_{k,k,t}^{i,j}\varphi_{b,k,t}^{i,j}  \quad \forall \; (i,j),\; l,\;t&&
\end{nospaceflalign} 
\begin{nospaceflalign}\nonumber
cc_{b,k,t}^{i,j} \leq \underline{cc}_{k,k,t}^{i,j}cc_{b,b,t}^{i,j} + \underline{cc}_{b,b,t}^{i,j}cc_{k,k,t}^{i,j} &&
\end{nospaceflalign} 
\begin{nospaceflalign}\nonumber
\quad - \underline{cc}_{b,b,t}^{i,j} \underline{cc}_{k,k,t}^{i,j}\varphi_{b,k,t}^{i,j}  \quad \forall \; (i,j),\; l,\;t&&
\end{nospaceflalign} 
\begin{nospaceflalign}\nonumber
cc_{b,k,t} = \sum\nolimits_{(i,j) \in \Omega^{cc}_{b,k,t}} cc_{b,k,t}^{i,j}  \quad \forall \; l,\; t &&
\end{nospaceflalign} 
\begin{nospaceflalign}
\underline{c}_{b,b,t}^{i} \varphi_{b,k,t}^{i,j} \leq c_{b,b,t}^{i,j} \leq \overline{c}_{b,b,t}^{i} \varphi_{b,k,t}^{i,j} \quad \forall \; (i,j),\; l,\; t &&
\end{nospaceflalign} 
\begin{nospaceflalign}\nonumber
c_{b,b,t} = \sum\nolimits_{(i,j) \in\Omega^{cc}_{b,k,t}} c_{b,b,t}^{i,j} \quad \forall \; l,\; t  &&
\end{nospaceflalign} 
\begin{nospaceflalign}\nonumber
\underline{c}_{k,k,t}^{i} \varphi_{b,k,t}^{i,j} \leq c_{k,k,t}^{i,j} \leq \overline{c}_{k,k,t}^{i} \varphi_{b,k,t}^{i,j} \quad \forall \; (i,j),\; l,\; t &&
\end{nospaceflalign} 
\begin{nospaceflalign}\nonumber
c_{k,k,t} = \sum\nolimits_{(i,j) \in\Omega^{cc}_{b,k,t}} c_{k,k,t}^{i,j} \quad \forall \; l,\; t  &&
\end{nospaceflalign} 
\begin{nospaceflalign}\nonumber
\sum\nolimits_{(i,j) \in \Omega^{cc}_{b,k,t}} \varphi_{b,k,t}^{i,j} =  1 \quad \forall \; l,\; t  &&
\end{nospaceflalign} 
\begin{nospaceflalign}\nonumber
\varphi_{b,k,t}^{i,j} \in \{0,1\} \quad \forall \; (i,j),\; l,\; t  &&
\end{nospaceflalign} 
\end{subequations}
where $(i,j)\in\Omega^{cc}_{b,k,t} \coloneqq [\underline{c}_{b,b,t}^{i}, \overline{c}_{b,b,t}^{i}] \times [\underline{c}_{k,k,t}^{j}, \overline{c}_{k,k,t}^{j}]$. Note that unique $c_{b,b,t}^{i,j}$ and $c_{k,k,t}^{i,j}$ variables must be introduced for every line $l$ where the under-estimators are constructed.

\subsubsection{Cycle Constraints}\label{sec:cycle_constraint}
In the second order cone relaxations used in $\mathsf{(M)}$ and $\mathsf{(Mf)}$, Kirchhoff's voltage law (KVL) is no longer guaranteed to be satisfied, but can be enforced through the \emph{cycle constraints},
\begin{nospaceflalign}\label{eq:SOCR_EQUIV_2}
\sum\nolimits_{(b,k)\in \mathcal{L}_c} \theta_{b,k,t} = 0 &&
 \end{nospaceflalign}
 for all $t$ and
\begin{nospaceflalign}\label{eq:SOCR_EQUIV_3}
 \theta_{b,k,t} = -\arctan(s_{b,k,t}/c_{b,k,t}) &&
 \end{nospaceflalign}
for all $l$ and $t$. In problem $\mathsf{(M)}$, these constraints are ignored (no refinement is necessary since solutions are enumerated with integer cuts). In problem $\mathsf{(Mf)}$, however, as the \emph{Inner Algorithm} iterates, these constraints are gradually enforced as needed by addition and refinement of piecewise outer approximations. We construct the respective piecewise under- and over-estimators for each $\theta_{b,k,t} = -\arctan (s_{b,k,t}/c_{b,k,t} ) $  term, where
\begin{subequations}
\renewcommand{\theequation}{$\mathsf{CC}$}
 \label{eq:CC}
\begin{nospaceflalign}\nonumber
\theta_{b,k,t}^{i,j} \geq \alpha^{i,j}_n s_{b,k,t}^{i,j} + \beta^{UE}_n c_{b,k,t}^{i,j} + \gamma^{UE}_n \quad\forall n,\; (i,j),\; l,\;t  &&
\end{nospaceflalign} 
\begin{nospaceflalign}\nonumber
\theta_{b,k,t}^{i,j}  \leq \alpha^{i,j}_n s_{b,k,t}^{i,j} + \beta^{OE}_n c_{b,k,t}^{i,j} + \gamma^{OE}_n \quad\forall n,\; (i,j),\; l,\;t  &&
\end{nospaceflalign} 
\begin{nospaceflalign}
\theta_{b,k,t} = \sum\nolimits_{(i,j) \in \Omega^{cs}_{b,k}} \theta_{b,k,t}^{i,j}&&
 \end{nospaceflalign}
 \begin{nospaceflalign}\nonumber
\sum\nolimits_{(b,k)\in \mathcal{L}_c} \theta_{b,k,t} = 0 &&
 \end{nospaceflalign}
\end{subequations}
where $n \in\{1,2\}$ and the parameters $\alpha, \beta,$ and $\gamma$ are based on the planes constructed in \cite{Kocuk2016}; then the bivariate piecewise partition is exact to the approach presented above in $(\mathsf{UE})$. Please see \cite{Liu2017a} for implementation details.

\subsubsection{Optimization-Based Bounds Tightening}\label{sec:obbt}
The optimization-based bounds tightening (OBBT) is only computed for the second-order cone variables $c_{b,k,t}$ and $s_{b,k,t}$ to perform domain reduction on the initial lower-bounding subproblem $(\mathsf{Mf})$. This approach results in two optimization routines per variable, i.e., 
\begin{subequations}
\begin{nospaceflalign}\nonumber
\underline{c}_{b,k,t} \leftarrow \max\big(\underline{c}_{b,k,t}, \min \{c_{b,k,t} | c(\mathsf{Mf}), z^0_U \leq z^*_U \}\big)  
\end{nospaceflalign} 
\begin{nospaceflalign}\nonumber
\overline{c}_{b,k,t} \leftarrow \min\big(\overline{c}_{b,k,t}, \max \{c_{b,k,t} | c(\mathsf{Mf}), z^0_U \leq z^*_U \}\big)
\end{nospaceflalign} 
\begin{nospaceflalign}\nonumber
\underline{s}_{b,k,t} \leftarrow \max\big(\underline{s}_{b,k,t}, \min \{s_{b,k,t} | c(\mathsf{Mf}), z^0_U \leq z^*_U \}\big)
\end{nospaceflalign} 
\begin{nospaceflalign}\nonumber
\overline{s}_{b,k,t} \leftarrow \min\big(\overline{s}_{b,k,t}, \max \{s_{b,k,t} | c(\mathsf{Mf}), z^0_U \leq z^*_U \}\big)  
\end{nospaceflalign} 
\end{subequations}
for all $l$ and $t$ where $c(\mathsf{Mf})$ denotes the constraint set of $(\mathsf{Mf})$. This procedure is computed selectively for $c_{b,k,t}$ and $s_{b,k,t}$ corresponding to large violations in second-order cone constraints \eqref{eq:SOCR_EQUIV_1}.
}

%
%
%
%
\section{Numerical Results}
\label{sec:results}

We now test our global \revision{UC-AC} solution algorithm on four benchmark problems: a $6$-bus test system (\texttt{6-bus}) with 3 generators \cite{Fu2006}, two $24$-bus test systems -- \texttt{RTS-79} \cite{subcommittee1979} and \texttt{RTS-96} -- each with 33 generators \cite{Wong1999}, and a modified IEEE $118$-bus test system (\texttt{IEEE-118mod}) with 54 generators \cite{Fu2006}. The scheduling horizon for all test cases is 24 hours at hourly time resolution. Our global solution algorithm is implemented in \emph{Pyomo}, a Python-based optimization modeling language \cite{Hart2017}. All computational experiments are conducted on a $64$-bit server comprising \revision{$24$} CPUs (Intel(R) Xeon(R) CPU E5-2697 v2 @ 2.70GHz) with $256$ GB of RAM. All SOCP and MISOCP subproblems are solved using Gurobi $6.5.2$ \cite{gurobi} limited to $24$ threads. All NLP subproblems are solved with Ipopt $3.12.6$ \cite{Wachter2006} using HSL's MA$27$ linear solver \cite{HSL2013}. 

\revision{In addition to having a tight and compact formulation for better performance in global solution frameworks,} convergence speed is also a function of other characteristics of the underlying numerical problem that impact computational difficulty, including formulation size and degeneracy / symmetry in the solution space. Typically, there is a large subset of solutions that are within an $\epsilon$-tolerance of an optimal-cost schedule. To balance computational burden with solution quality, we initially set the Gurobi MIP gap to $0.1\%$. Then, if the optimality gap of our global solution algorithm does not show improvements within $N$ iterations, we tighten the MIP gap by a factor of $10$. For the nested algorithm for subproblem, the $\epsilon$ is also set to $0.1\%$.

In all of our computational experiments, we set $N=5$ with a total wall clock time limit of $14400s$ and a major iteration limit $q=30$. The optimality tolerance for both our global solution algorithm and its nested multi-tree algorithm are set to $0.1\%$.




\subsection{Computational Performance}

\begin{table*}[h]
  \begin{center}
    \caption{Numerical results for our global UC-AC solution algorithm}
    \label{table:UC_Results}
    \begin{tabular}{cccccc}
    \toprule
    Case & Upper Bound ($\$$) & Lower Bound ($\$$) & Optimality Gap ($\%$) & Wall Clock Time (s) & Iteration ($k$)\\    
    \midrule
    \texttt{6-bus} &  $101,763$ & $101,740$ & $0.02\%$ & $8.5$ & $2$\\
    \texttt{RTS-79} &  $895,040$ & $894,392$ & $0.07\%$ & $1394$ & $6$\\
    \texttt{RTS-96} &  $886,362$ & $885,707$ & $0.07\%$ & $321.0$ & $1$\\
    \texttt{IEEE-118mod} &  $835,926$ & $833,057$ & $0.34\%$ & $14400^*$ & $2$ \\
    \bottomrule
    \end{tabular}
  \end{center}
\end{table*}

Computational results for our global solution algorithm on the 4 benchmark problems are reported in Table \ref{table:UC_Results}. \revision{The second column reports the best obtained upper bound, which corresponds to the \textit{best known solution} to the \revision{UC-AC} problem. The third column reports the best obtained lower bound, which corresponds to the solution of the problem defined in $(\mathsf{M})$.} The relative optimality gap is shown in the fourth column, followed by the total wall clock time and the number of major iterations. All problems are solved to within a $0.5\%$ global optimality gap in under the wall clock time limit. For \texttt{IEEE-118mod}, we obtained a $0.34\%$ optimality gap after the first iteration (in approximately $8400 s$), which remains unchanged before the time limit is reached in major iteration $k=2$ with a $0.11\%$ MIP gap for the lower-bounding problem.

\revision{We also attempted to solve these UC-AC problems using the version 16.12.7 of the commercially available general MINLP solver, BARON\cite{Sahinidis1996,baron2005}. This general algorithm was not able to solve any of the UC-AC problems within a time limit of 10 hours. For the 6-bus case study, no significant progress was made in either the upper or lower bound with 13797 iterations of the algorithm. For these tests, we used default values for all algorithm tuning parameters.  For subproblem solvers, CLP/CBC was used for LP and MIP problems while IPOPT and FILTERSD were used for the NLP subproblems. It is possible that better performance could be obtained by additional tuning.}
  
  We also note that in contrast to research on global solution of MIP models, in which accepted optimality tolerances are typically $1\cdot10^{-4}$, standards for global solution of MINLP models are typically within $1$\% -- due to the relative increase in computational difficulty and maturity of global NLP subproblem solvers.

\subsection{Globally Optimal Unit Commitment Schedules}

Globally optimal schedules for our test cases are shown in Table \ref{table:Commitment_Case6}, \ref{table:Commitment_Case24}, and \ref{table:Commitment_Case118}; there are multiple globally optimal solutions for \texttt{RTS-79}, \texttt{RTS-96}, and \texttt{IEEE-118mod} (not reported here). The multiple global solutions are due to the symmetry, e.g., co-location of identical generating units at a given bus in the $24$-bus case and identical branches in the $118$-bus case. To partially remedy this problem, symmetry-breaking methods, e.g. see \cite{Ostrowski2012b},  can be incorporated within the proposed global solution algorithm for the UC-AC formulation. Relative to our local method for the unit commitment with AC transmission constraints \cite{castillo2016}, we observe that our global solution algorithm locate the same solution to \texttt{6-Bus}, a slightly improved solution to \texttt{RTS-79}, and a significantly improved solution to \texttt{IEEE-118mod}.

\begin{table}
  \begin{center}
    \caption{Commitments for the \texttt{6-Bus} System}\label{table:Commitment_Case6}
    \begin{tabular}{cc|c}
    \toprule
   Bus & Gen& Commitment (h) \\
    \hline
   B1 & G1 & $1$-$24$ \\
   B2 & G2 & $1$, $12$-$21$ \\
   B6 & G3 & $10$-$22$ \\
    \bottomrule
    \end{tabular}
  \end{center}
\end{table}

\begin{table}
  \begin{center}
    \caption{Commitments for the 24-Bus Systems}\label{table:Commitment_Case24}
    \begin{tabular}{cc|cc}
    \toprule
   Bus & Gen & \multicolumn{2}{c}{Commitment (h)} \\
   & & \texttt{RTS-79} & \texttt{RTS-96} \\
    \hline
   B1 & G1, G2 & $\O$ & $\O$ \\
   B1 & G3, G4 & $8$-$23$ & $8$-$23$ \\
   B2 & G5, G6& $10$ & $\O$ \\
   B2 & G7 & $8$-$24$ & $8$-$24$ \\
   B2 & G8 & $8$-$23$ & $8$-$23$ \\
   B7 & G9 & $1$-$23$ & $1$-$23$ \\
   B7 & G10 & $9$-$24$ & $10$-$24$ \\ 
   B7 & G11 & $10$-$18$ & $\O$ \\
   B13 & G12 & $11$-$22$ & $1$-$18$ \\
   B13 & G13 & $\O$ & $11$-$22$ \\
   B13 & G14 & $\O$ & $\O$ \\
   B14 & G15 & $1$-$24$ & $1$-$24$ \\   
   B15 & G16-G18  & $10$-$15$ & $\O$ \\
   B15 & G19, G20 & $10$-$13$ & $\O$ \\      
   B15 & G21 & $9$-$24$ & $9$-$24$ \\        
   B16 & G22 & $1$-$24$ & $1$-$24$ \\      
   B18 & G23 & $1$-$24$ & $1$-$24$ \\          
   B21 & G24 & $1$-$24$ & $1$-$24$ \\     
   B22 & G25-G30 & $1$-$24$ & $1$-$24$ \\     
   B23 & G31-G33 & $1$-$24$ & $1$-$24$ \\     
    \bottomrule
    \end{tabular}
  \end{center}
\end{table}

\begin{table}
  \begin{center}
    \caption{Commitments for the \texttt{IEEE-118mod} System}\label{table:Commitment_Case118}
    \begin{tabular}{cc|cc}
    \toprule
    Gen & Commitment (h)& Gen & Commitment (h) \\
    \hline  
   G1 & $\O$ & G28 & $1$-$24$ \\        
   G2 & $\O$ & G29 & $1$-$24$ \\
   G3 & $\O$ & G30 & $1$-$24$  \\
   G4 & $1$-$10$, $24$ & G31 & $\O$\\
   G5 & $1$-$24$ & G32 & $\O$ \\
   G6 & $\O$ & G33 & $\O$ \\
   G7 & $11$-$22$ & G34 & $7$-$24$ \\     
   G8 & $\O$ & G35 & $1$-$24$ \\                                                                                                                                                                                                                                                                                                   
   G9 & $\O$ & G36 & $1$-$24$ \\                                                                                                                                                                                                                                                                                                    
   G10 & $1$-$2$, $12$-$24$ & G37 &$8$-$23$ \\
   G11 & $1$-$24$ & G38 & $\O$ \\                                                                                                                                                                                                                                                                                                    
   G12 & $\O$ & G39 & $\O$ \\                                                                                                                                                                                                                                                                                                    
   G13 & $\O$ & G40 & $1$-$10$, $22$-$24$ \\ 
   G14 & $10$-$22$ & G41 & $\O$ \\      
   G15 & $\O$ & G42 & $\O$ \\                                                                                                                                                                                                                                                                                                   
   G16 & $9$-$16$ & G43 & $1$-$24$ \\
   G17 & $\O$ & G44 & $\O$ \\                                                                                                                                                                                                                                                                                                                                                                                                                                                                                                                                                                                                                                                     
   G18 & $\O$ & G45 & $1$-$24$ \\   
   G19 & $\O$ & G46 & $\O$ \\   
   G20 & $1$-$24$ & G47 & $\O$ \\                                                                                                                                                                                                                                                                                                                                                                                                                                                                                                                                                                                                                                                                  
   G21 & $8$-$24$ & G48 & $\O$ \\                                                                                                                                                                                                                                                                                                                                                                                                                                                                                                                                                                                                                                                     
   G22 & $\O$ & G49 & $\O$ \\                                                                                                                                                                                                                                                                                                         
   G23 & $\O$ & G50 & $\O$ \\                                                                                                                                                                                                                                                                                                 
   G24 & $9$-$23$ & G51 & $9$-$13$ \\        
   G25 & $\O$ & G52 & $14$-$23$ \\  
   G26 & $\O$ & G53 & $7$-$24$ \\                                                                                                                                                                                                                                                                                                                                                                                                                                                                                                                                                                                                                                                                  
   G27 & $1$-$2$, $13$-$24$ & G54 & $9$-$23$ \\                                                                                                                                                                                                                                                                                                                                                                                                                                                                                                                                                                                                                                                         
    \bottomrule
    \end{tabular}
  \end{center}
\end{table}

\section{Conclusions} 
\label{sec:conclusions}
\revision{Solving the UC-AC problem is fundamental to solving real-world operations and market settlements that fully incorporate the impact of alternating current physics on the network.} We have introduced, to the best of our knowledge, the first such approach to solving this practically critical and computationally difficult problem. Although our obtained run times are still longer than those required for operations, our proposed approach can be used to assess the provably (near-) global optimality of ``off-line" solutions as well as test and validate other algorithmic approaches including heuristics and local solution techniques, e.g. see \cite{castillo2016}.

Future directions for research include \revision{improving relaxations of the UC-AC}, incorporating symmetry-breaking methods, \revision{and other cutting plane techniques to improve the efficiency in solving the mixed-integer master problem; improvements to the mixed-integer refinement problem in the nested algorithm include adaptive, non-uniform partitioning schemes. Security considerations and uncertainties do not alter the core UC-AC problem that needs to be solved, but does increase the dimensionality of the problem; such dimensionality increase is addressable through decomposition and parallelization techniques, which are extensions left for future work.}

\section*{Acknowledgments}
Sandia National Laboratories is a multimission laboratory managed and operated by National Technology \& Engineering Solutions of Sandia, LLC, a wholly owned subsidiary of Honeywell International Inc., for the U.S. Department of Energy’s National Nuclear Security Administration under contract DE-NA0003525.

\emph{Disclaimer:} This paper describes objective technical results and analysis. Any subjective views or opinions that might be expressed in the paper do not necessarily represent the views of the U.S. Department of Energy or the United States Government.


\bibliography{jianfeng_optimization,jianfeng_uc_acopf}                                   
\newpage

\include{NCUC_Appendix}

\end{document}